\documentclass{agtart_a}
\pdfoutput=1

\usepackage{pinlabel}

\title{Ordering the Reidemeister moves of a classical knot}

\author{Alexander Coward}
\givenname{Alexander}
\surname{Coward}
\address{Mathematical Institute\\24-29 St Giles'\\Oxford\\OX1 3LB\\UK}
\email{cowardaw@maths.ox.ac.uk}
\urladdr{}

\volumenumber{6}
\issuenumber{}
\publicationyear{2006}
\papernumber{25}
\startpage{659}
\endpage{671}

\doi{}
\MR{}
\Zbl{}


\received{8 June 2005}
\revised{13 April 2006}
\accepted{17 April 2006}
\published{18 May 2006}
\publishedonline{18 May 2006}
\proposed{}
\seconded{}
\corresponding{}
\editor{}
\version{}

\arxivreference{}  

\AtBeginDocument{\let\bar\wbar\let\tilde\wtilde\let\hat\what}
\def\circleright{\hspace{7pt}\raise 3pt\hbox{\circle{8}}\hspace{-4pt}{\to}\,}

\newtheorem{thm}{Theorem}
\newtheorem{prop}{Proposition}
\newtheorem{lem}{Lemma}
\newtheorem{cltwo}{Corollary to Lemma 2}
\makeautorefname{cltwo}{Corollary to Lemma 2}

\newtheorem{clthree}{Corollary to Lemma 3}
\makeautorefname{clthree}{Corollary to Lemma 3}

\theoremstyle{definition} 
\newtheorem{dfn}{Definition}

\begin{document}

\begin{abstract}
We show that any two diagrams of the same knot or link are connected
by a sequence of Reidemeister moves which are sorted by type.
\end{abstract}

\maketitle

It is one of the founding theorems of knot theory that any two
diagrams of a given link may be changed from one into the other by a
sequence of Reidemeister moves. One of the reasons why this result
is so crucial to the subject is that it allows one to define a link
invariant as an invariant of a diagram which is unchanged under
Reidemeister moves.

\begin{figure}[ht!]\small\hair 2pt
\labellist
\pinlabel {$\Omega_1^\uparrow$} [b] <2pt,0pt> at 50 21
\pinlabel {$\Omega_1^\downarrow$} [t] <2pt,0pt> at 50 14
\pinlabel {$\Omega_2^\uparrow$} [b] <2pt,0pt> at 164 21
\pinlabel {$\Omega_2^\downarrow$} [t] <2pt,0pt> at 164 14
\pinlabel {$\Omega_3$} [b] <1pt,0pt> at 279 19
\endlabellist
\centering
\includegraphics{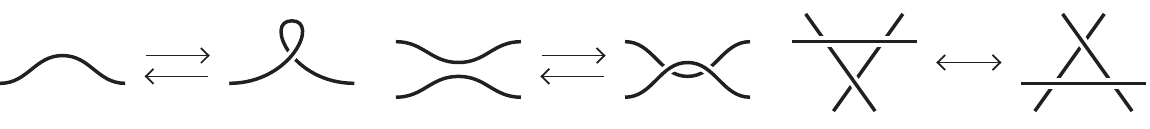}
\caption{Reidemeister moves} \label{F:1}
\end{figure}
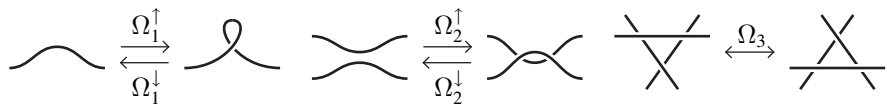

Since Reidemeister's seminal paper on this topic in 1927 \cite{RM},
there have been a number of steps taken to strengthen the original
result in a variety of directions. In 1983, Bruce Trace \cite{Tr}
proved that type 1 moves may be omitted in the case where two knot
diagrams have the same winding number and framing. Recent work by
Joel Hass and Jeffrey Lagarias \cite{HL} has placed a bound on the
number of moves required when one of the diagrams is the trivial
unknot diagram.

In this paper we shall address the question of whether, given any
two diagrams of a knot or link, there exists a sequence of
Reidemeister moves between them which is sorted by type. We answer
this question in the affirmative with the following theorem:

\begin{thm}\label{thm:1}
Given two diagrams $D_1$ and $D_2$ for a
link $L$, $D_1$ may be turned into $D_2$ by a sequence of
$\Omega_1^\uparrow$ moves, followed by a sequence of
$\Omega_2^\uparrow$ moves, followed by a sequence of $\Omega_3$
moves, followed by sequence of $\Omega_2^\downarrow$ moves.

Furthermore, if $D_1$ and $D_2$ are diagrams of a link where the
winding number and framing of each component is the same in each
diagram, then $D_1$ may be turned into $D_2$ by a sequence of
$\Omega_2^\uparrow$ moves, followed by a sequence of $\Omega_3$
moves, followed by a sequence of $\Omega_2^\downarrow$ moves.
\end{thm}

In this paper, all link diagrams shall be regarded as 4--valent
graphs embedded in $\mathbb{R}^2$ with signed intersections to
denote overcrossings or undercrossings. All diagrams shall be
oriented so as to represent an oriented link. $\Omega_1^\uparrow$,
$\Omega_1^\downarrow$, $\Omega_2^\uparrow$, $\Omega_2^\downarrow$
and $\Omega_3$ shall denote Reidemeister moves where the arrow
indicates whether the move increases the number of crossings in the
diagram, or decreases it, as shown in \fullref{F:1}. The winding number
of a component of a link in a diagram is intuitively speaking the
number of times that one must rotate anticlockwise when walking once
around that component in the specified orientation. The framing
(also known as the writhe) of a knot diagram is the number of
crossings where the upper strand's orientation is 90 degrees
clockwise from that of the lower strand, minus the number of
crossings where the upper strand's orientation is 90 degrees
anticlockwise from that of the lower strand. The framing of a
component of a link diagram is obtained by taking the difference
over only those crossings where both strands belong to the component
in question. For more on these notions see Trace \cite{Tr}.

Returning to \fullref{thm:1}, the first part of the theorem in fact
follows from the second part because of the following proposition:

\begin{prop}\label{prop:1}
Let $D_1$ and $D_2$ be two diagrams
for a link $L$. Then we may apply $\Omega_1^\uparrow$ moves to $D_1$
so as to obtain a new diagram $D_1'$ with the all the same winding
numbers and framings as $D_2$.
\end{prop}

\proof We know that $D_1$ may be changed into $D_2$ by a
sequence of Reidemeister moves. Note that only $\Omega_1$ moves
change the winding numbers and framings of a diagram, and that each
$\Omega_1$ move changes the winding number of the component on which
it acts by $\pm1$ and the framing of the component on which it acts
by $\pm1$. For each $\Omega_1$ move in a sequence of Reidemeister
moves from $D_1$ and $D_2$, we may carry out a $\Omega_1^\uparrow$
move on an edge in the same component in $D_1$ with the same effect
on the winding number and framing of that component. After
completing these moves we will have a new diagram $D_1'$ with the
same winding numbers and framings as $D_2$. \endproof

We shall now turn our attention to the proof of \fullref{thm:1}. Our
strategy will be to simulate each $\Omega_3$ move with a sequence of
$\Omega_2^\uparrow$ moves. In order to achieve this we will need to
develop some new notation.

\begin{figure}[ht!]\small\hair 7pt
\labellist
\pinlabel {$D\cup C$} [t] <0pt,0pt> at 67 20
\pinlabel {$D'$} [t] <0pt,0pt> at 237 19
\endlabellist
\centering
\includegraphics{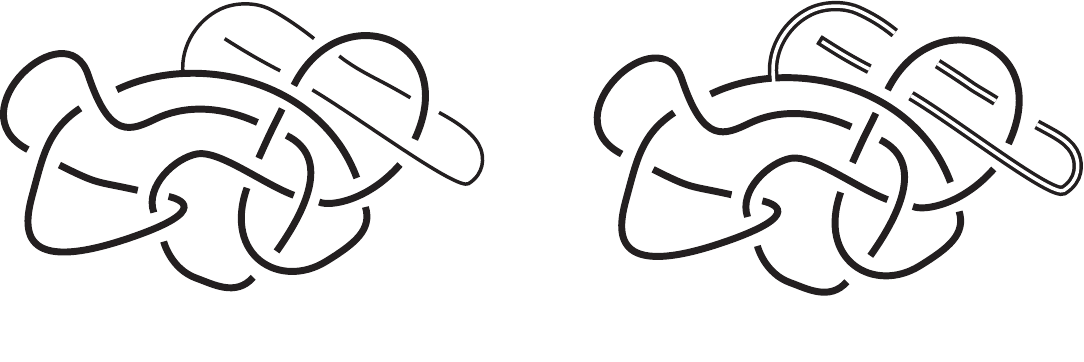}
\caption{Adding a tail} \label{F:2}
\end{figure}

\begin{dfn}\label{def:1}
Let $D$ be a link diagram in $\mathbb{R}^2$
and let $c\co I\rightarrow \mathbb{R}^2$ be an embedded path such that
$c(0)$ lies on the interior of an edge of $D$, $c(1)$ does not lie
on $D$, and where $c(\textrm{int}(I))$ intersects $D$ transversely
in a finite number of points none of which are vertices of $D$, and
which are given signings to denote whether the path crosses above or
below $D$. Let $C$ denote the image of this path so that $D \cup C$
is a graph which is 3--valent at one vertex, 1--valent at another, and
4--valent otherwise as shown in the left hand image of \fullref{F:2}.
\end{dfn}

Let $C \times [-\epsilon,\epsilon]$ denote a small product
neighbourhood of $C$ such that $(C \times [-\epsilon,\epsilon]) \cap
D = (C \cap D) \times [-\epsilon,\epsilon]$. Let $D'$ be the link
diagram whose 4--valent graph is $$D \cup \partial(C \times
[-\epsilon,\epsilon])\backslash(c(0) \times (-\epsilon,\epsilon))$$
and whose vertex signings are induced by those of $D \cup C$. The
orientation of $D'$ shall be that induced by $D$.

We shall say that \emph{$D'$ is obtained from $D$ by adding a tail
along $C$}. We shall call $\partial(C \times
[-\epsilon,\epsilon])\backslash(c(0) \times (-\epsilon,\epsilon))$
\emph{the tail in $D'$} and we shall refer to \emph{$C$ as the core
of this tail}. We shall write $D\rightsquigarrow{}D'$. Note that if
$D\rightsquigarrow{}D'$ then $D'$ may be obtained from $D$ by a
sequence of $\Omega_2^\uparrow$ moves. Note also that the core of a
tail is an embedded arc, and not an immersed one.

\begin{figure}[ht!]\small\hair 7pt
\labellist
\pinlabel {$D_2$} [t] <0pt,0pt> at 46 18
\pinlabel {$D_3$} [t] <0pt,0pt> at 226 18
\hair 2pt
\pinlabel {$\Omega_2^\uparrow$\, $\Omega_2^\uparrow$\, $\Omega_3$} [bl] <-5pt,0pt> at 112 54
\endlabellist
\centering
\includegraphics{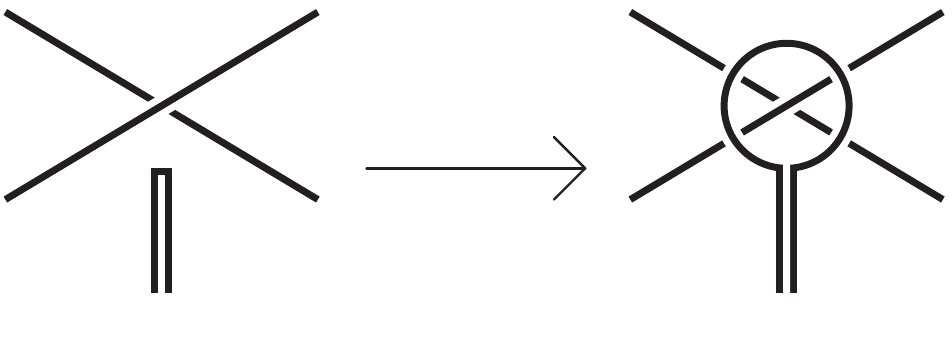}
\caption{Turning a tail into a lollipop} \label{F:3}
\end{figure}

\begin{dfn}\label{def:2}
Suppose that $D_1\rightsquigarrow{}D_2$. Let
$D_3$ be obtained from $D_2$ by performing two $\Omega_2^\uparrow$
moves and one $\Omega_3$ move as shown in \fullref{F:3}. We shall then
say that \emph{$D_3$ is obtained from $D_1$ by adding a lollipop}
and we shall write $D_1 \circleright D_3$.\end{dfn}

Later on it will be important to distinguish between the part of the
lollipop which circles the crossing and the part which consists of
two parallel strands. We shall call these the \emph{circle part} and
the \emph{tail part} of the lollipop respectively.

We are now in a position to say how we are going to simulate
$\Omega_3$ moves by means of $\Omega_2^\uparrow$ moves. This is
captured in the following important lemma.

\begin{lem}\label{lem:1}
Suppose that $D_2$ is obtained from $D_1$ by
means of an $\Omega_3$ move. Then we may construct a diagram $D_3$
such that:\begin{enumerate}
    \item $D_3$ may be obtained from $D_1$ by a sequence of $\Omega_2^\uparrow$ moves.
    \item $D_2\circleright{}D_3$
\end{enumerate}
\end{lem}

\proof Let $D_3$ be as shown below.
$$\labellist\small
\pinlabel {$D_1$} <0pt,2pt> at 56 25
\pinlabel {$D_2$} <2pt,2pt> at 128 25
\pinlabel {$D_3$} <0pt,2pt> at 91 62
\pinlabel {seq$(\Omega_2^\uparrow)$} [br] <0pt,0pt> at 71 44
\endlabellist
\includegraphics{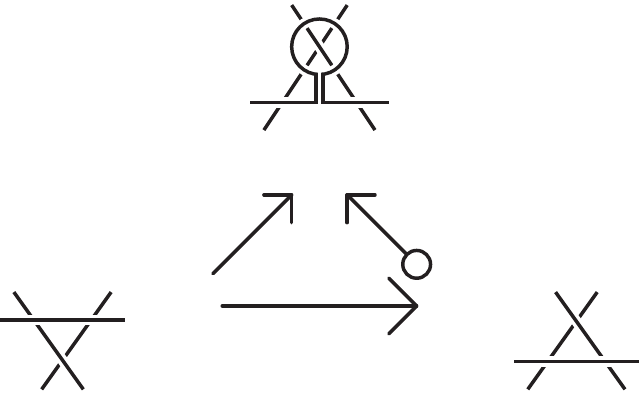}\proved$$

\begin{lem}\label{lem:2}
Let $D_1\rightsquigarrow{}D_1'$. Suppose that
$D_2$ may be obtained from $D_1$ by a Reidemeister move of type
$\Omega_2^\uparrow$. It is possible to construct a diagram $D_2'$
such that:\begin{enumerate}
    \item $D_2'$ may be obtained from $D_1'$ by a sequence of $\Omega_2^\uparrow$ moves.
    \item $D_2\rightsquigarrow{}D_2'$
\end{enumerate}\end{lem}

\proof Let $C$ denote the core of the tail in $D_1'$. Let
$E_1$ and $E_2$ be the edges (not necessarily distinct) in $D_1$
upon which our $\Omega_2^\uparrow$ move takes place. Note that $E_1$
and $E_2$ are incident to a face $F$ of the diagram $D_1$. Let $x_1$
(resp. $x_2$) be a point on $E_1$ (resp. $E_2$) which does not lie
in $C \times [-\epsilon,\epsilon]$. Let $P$ be an embedded path from
$x_1$ to $x_2$ whose interior lies entirely in $F$ and which crosses
$C \times [-\epsilon,\epsilon]$ transversely in a finite number of
intervals. Let $P'$ be a path obtained from $P$ by extending it a
small amount at $x_2$ into the neighbouring face. Then $D_2'$ may be
formed by adding a tail along $P'$ to $D_1'$ as shown below.
$$\labellist\small\hair 5pt
\pinlabel {$D_1'$} [r] <0pt,0pt> at 107 76
\pinlabel {$D_2'$} [l] <0pt,0pt> at 153 76
\pinlabel {$D_1$} [r] <0pt,0pt> at 107 39
\pinlabel {$D_2$} [l] <0pt,0pt> at 153 39
\pinlabel \scalebox{2}{\rotatebox{90}{$\rightsquigarrow$}} [r] <-2pt,2pt> at 107 57
\pinlabel \scalebox{2}{\rotatebox{90}{$\rightsquigarrow$}} [l] <0pt,2pt> at 153 57
\hair 3pt
\pinlabel {seq$(\Omega_2^\uparrow)$} [b] <0pt,0pt> at 126 76
\pinlabel {$\Omega_2^\uparrow$} [t] <0pt,0pt> at 126 39
\endlabellist
\includegraphics{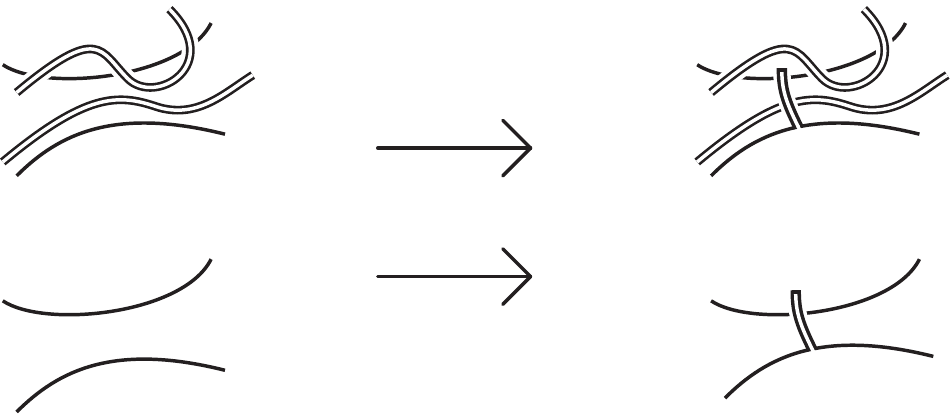}\proved$$

\begin{cltwo}\label{cor:2}
Let $D_1\rightsquigarrow{}D_1'$.
Suppose that $D_2$ may be obtained from $D_1$ by means of a sequence
of Reidemeister moves of type $\Omega_2^\uparrow$. It is possible to
construct a diagram $D_2'$ such that:\begin{enumerate}
    \item $D_2'$ may be obtained from $D_1'$ by a sequence of $\Omega_2^\uparrow$ moves.
    \item $D_2\rightsquigarrow{}D_2'$
\end{enumerate}\end{cltwo}

\proof Let $D_1 = E_1, \ldots , E_n = D_2$ be a sequence
of diagrams such that
$E_i\stackrel{\Omega_2^\uparrow}{\longrightarrow}E_{i+1}$. Thus we
have:
$$\labellist\small\hair 5pt
\pinlabel {$D_1'=E_1'$} [r] <0pt,0pt> at 47 55
\pinlabel {$D_1=E_1$} [r] <0pt,0pt> at 47 15
\pinlabel {$E_2'$} <0pt,0pt> at 107 55
\pinlabel {$E_2$} <0pt,0pt> at 107 15
\pinlabel {$E_n=D_2$} [l] <0pt,0pt> at 236 15
\pinlabel {$E_n'=D_2'$} [l] <0pt,0pt> at 236 55
\pinlabel \scalebox{2}{\rotatebox{90}{$\rightsquigarrow$}} [r] <-2pt,2pt> at 47 35
\pinlabel \scalebox{2}{\rotatebox{90}{$\rightsquigarrow$}}  <0pt,0pt> at 107 35
\pinlabel \scalebox{2}{\rotatebox{90}{$\rightsquigarrow$}} [l] <0pt,2pt> at 236 35
\hair 3pt
\pinlabel {seq$(\Omega_2^\uparrow)$} [b] <0pt,0pt> at 65 55
\pinlabel {seq$(\Omega_2^\uparrow)$} [b] <-2pt,0pt> at 143 55
\pinlabel {seq$(\Omega_2^\uparrow)$} [b] <-2pt,0pt> at 210 55
\pinlabel {$\Omega_2^\uparrow$} [t] <0pt,0pt> at 65 14
\pinlabel {$\Omega_2^\uparrow$} [t] <0pt,0pt> at 143 14
\pinlabel {$\Omega_2^\uparrow$} [t] <0pt,0pt> at 210 14
\pinlabel {$\cdots$} <0pt,0pt> at 180 14
\pinlabel {$\cdots$} <0pt,0pt> at 180 55
\endlabellist
\includegraphics{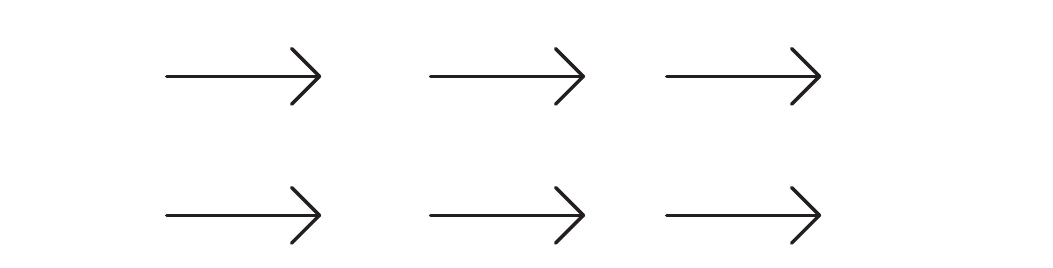}\proved$$

There is a similar pair of results for the adding of lollipops:

\begin{lem}\label{lem:3}
Let $D_1\circleright{}D_1'$. Suppose that $D_2$
may be obtained from $D_1$ by a Reidemeister move of type
$\Omega_2^\uparrow$. It is possible to construct a diagram $D_2'$
such that:\begin{enumerate}
    \item $D_2'$ may be obtained from $D_1'$ by a sequence of $\Omega_2^\uparrow$ moves.
    \item $D_2\circleright{}D_2'$
\end{enumerate}\end{lem}

\proof In this case we proceed exactly as in the proof of
\fullref{lem:1} except that we insist that the path $P$ avoids the circle
part of of the lollipop.
$$\labellist\small\hair 5pt
\pinlabel {$D_1'$} <2pt,0pt> at 3 55
\pinlabel {$D_1$} <2pt,0pt> at 3 15
\pinlabel {$D_2$} <0pt,0pt> at 80 15
\pinlabel {$D_2'$} <0pt,0pt> at 80 55
\hair 3pt
\pinlabel {seq$(\Omega_2^\uparrow)$} [b] <0pt,0pt> at 38 55
\pinlabel {$\Omega_2^\uparrow$} [t] <0pt,2pt> at 38 15
\endlabellist
\includegraphics{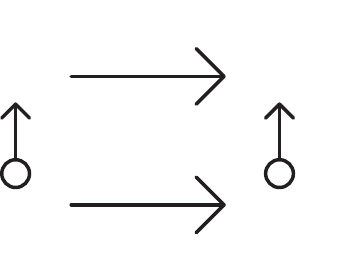}\proved$$

\begin{clthree}\label{cor:3}
Let $D_1\circleright{}D_1'$. Suppose
that $D_2$ may be obtained from $D_1$ by sequence of Reidemeister
moves of type $\Omega_2^\uparrow$. It is possible to construct a
diagram $D_2'$ such that:\begin{enumerate}
    \item $D_2'$ may be obtained from $D_1'$ by a sequence of $\Omega_2^\uparrow$ moves.
    \item $D_2\circleright{}D_2'$
\end{enumerate}\end{clthree}

\proof As before let $D_1 = E_1, \ldots , E_n = D_2$ be a
sequence of diagrams such that
$E_i\stackrel{\Omega_2^\uparrow}{\longrightarrow}E_{i+1}$. In this
case we have:
$$\labellist\small\hair 5pt
\pinlabel {$D_1'=E_1'$} [r] <0pt,0pt> at 47 55
\pinlabel {$D_1=E_1$} [r] <0pt,0pt> at 47 15
\pinlabel {$E_2'$} <0pt,0pt> at 107 55
\pinlabel {$E_2$} <0pt,0pt> at 107 15
\pinlabel {$E_n=D_2$} [l] <0pt,0pt> at 236 15
\pinlabel {$E_n'=D_2'$} [l] <0pt,0pt> at 236 55
\hair 3pt
\pinlabel {seq$(\Omega_2^\uparrow)$} [b] <0pt,0pt> at 65 55
\pinlabel {seq$(\Omega_2^\uparrow)$} [b] <-2pt,0pt> at 143 55
\pinlabel {seq$(\Omega_2^\uparrow)$} [b] <-2pt,0pt> at 210 55
\pinlabel {$\Omega_2^\uparrow$} [t] <0pt,0pt> at 65 14
\pinlabel {$\Omega_2^\uparrow$} [t] <0pt,0pt> at 143 14
\pinlabel {$\Omega_2^\uparrow$} [t] <0pt,0pt> at 210 14
\pinlabel {$\cdots$} <0pt,0pt> at 180 14
\pinlabel {$\cdots$} <0pt,0pt> at 180 55
\endlabellist
\includegraphics{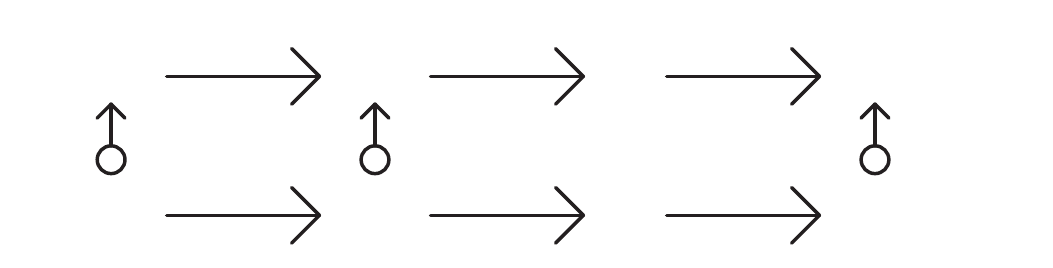}\proved$$

We need one more result before we can turn to the proof of \fullref{thm:1}.

\begin{prop}\label{prop:2}
Suppose that $D_2$ is obtained from $D_1$
by the addition of a sequence of tails and lollipops. Then there
exists a diagram $D_3$ such that:\begin{enumerate}
    \item $D_3$ may be obtained from $D_2$ by means of a sequence of
    $\Omega_2^\uparrow$ moves followed by a sequence of $\Omega_3$
    moves.
    \item $D_3$ may be obtained from $D_1$ by means of a sequence of
    $\Omega_2^\uparrow$ moves.
\end{enumerate}\end{prop}

\proof Our strategy will be to construct $D_3$ from $D_2$
in accordance with the first condition and then show that our new
diagram $D_3$ satisfies the second condition.

\begin{figure}[ht!] \centering\small
\labellist
\pinlabel $P$ [b] at 26 112
\pinlabel $x$ [t] at 82 102
\pinlabel $v$ [b] at 104 99
\endlabellist
\includegraphics{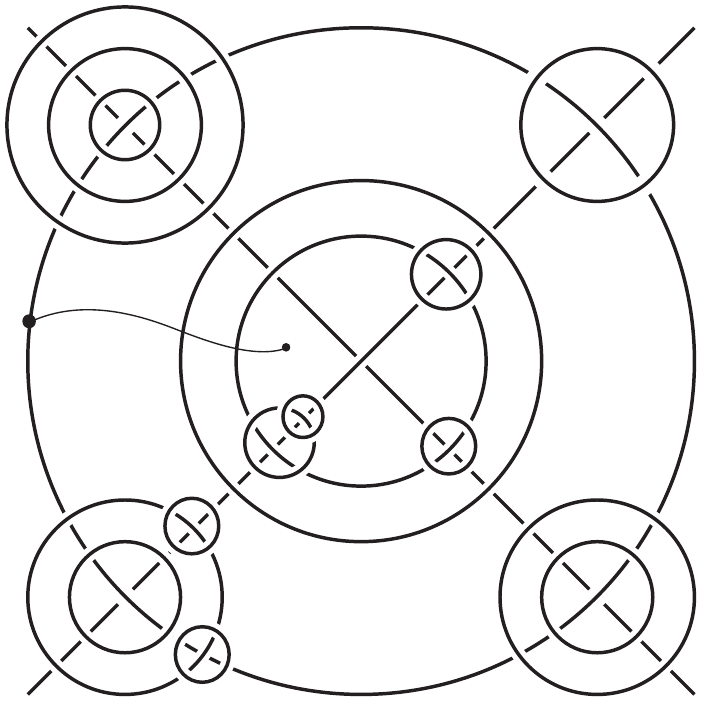}
\caption{Diagram showing the circle part of a lollipop in $D_2$}\label{F:4}
\end{figure}

Let $D_1 = E_1 , \ldots , E_n = D_2$ be a sequence of diagrams such
that for each $i \in \{1 , \ldots , n-1 \}$ either $E_i
\rightsquigarrow E_{i+1}$ or $E_i \circleright E_{i+1}$. We shall
start by performing moves on the circle parts of the lollipops in
$D_2$. Note that each of these circle parts is associated with a
particular vertex of $D_2$, namely the vertex around which the
circle part was originally added, and furthermore that the circle
parts associated to a particular vertex are disjoint and concentric.
Consider all the circle parts around a vertex $v$. Let $x$ be some
point in a region $R$ of $D_2$ which neighbours $v$. Let $P$ be a
path from $x$ to a point on the outermost circle part $C$ associated
with $v$ which avoids circle parts of other lollipops and avoids the
tail part of $C$, as shown in \fullref{F:4} which omits any tail parts of
$D_2$ for the sake of clarity.

\begin{figure}[ht!] \centering
\includegraphics{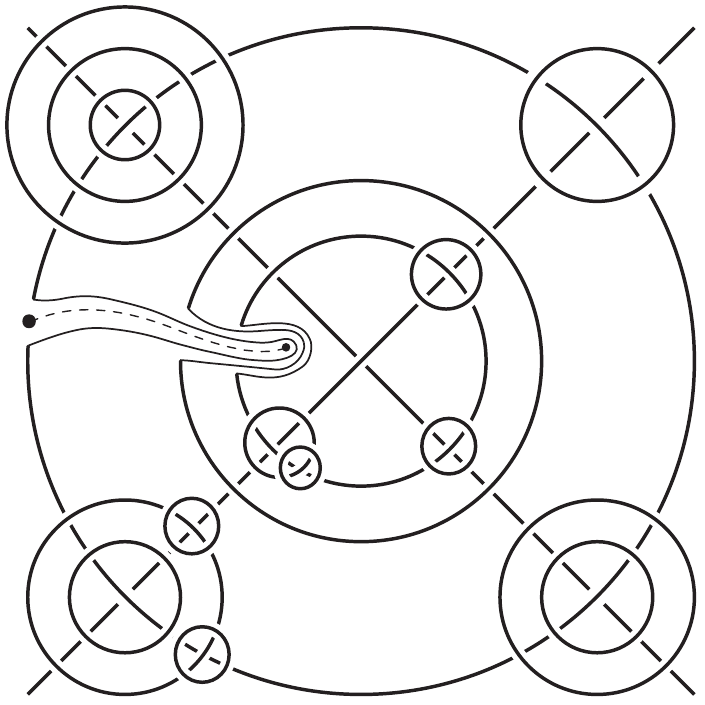}
\caption{Add tails along $P$ using $\Omega_2^\uparrow$ moves}\label{F:5}
\end{figure}

We may now undertake a sequence of type $\Omega_2^\uparrow$ moves in
a neighbourhood of $P$ as follows. It will be convenient to use the
language of adding tails, but one should view this as a shorthand
for describing a sequence of $\Omega_2^\uparrow$ moves. First add a
tail to the innermost circle part associated to $v$ along the part
of $P$ which is inside that circle part. Note that $P$ will be
disjoint from this tail in the resulting diagram apart from at $x$.
Extend the tail a small amount so that $P$ and the tail are now
disjoint. Continue by adding tails to all the circle parts
associated to $v$ along $P$ in the same way, working in order from
the innermost circle part to the outermost circle part, $C$. This
procedure is illustrated in \fullref{F:5}.

It is worth remembering the tail part of the diagram not shown in
the figure, and observing that as long as we are just performing
$\Omega_2^\uparrow$ moves then we can simply go over that part as
required.

\begin{figure}[ht!] \centering
\includegraphics{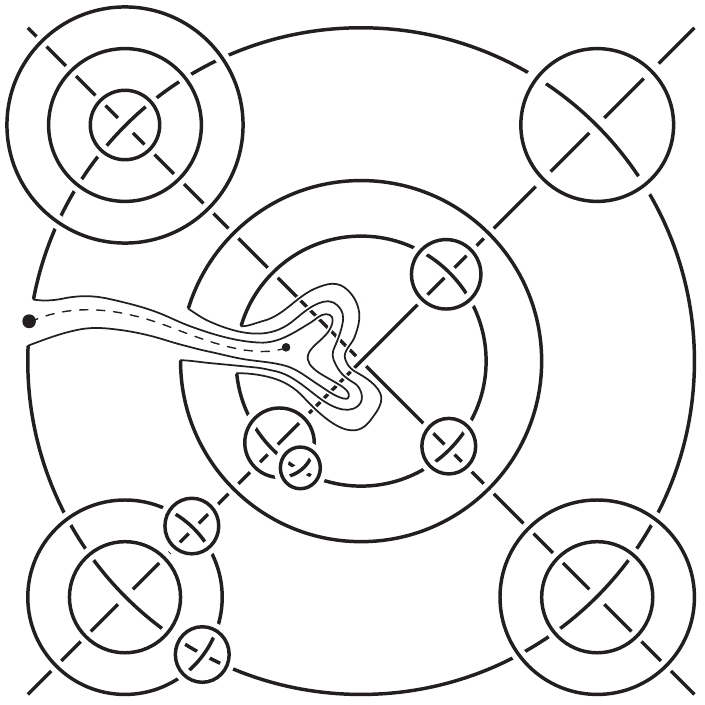}
\caption{Perform some more $\Omega_2^\uparrow$ moves near each
vertex $v$}\label{F:6}
\end{figure}

Since $x$ was chosen to lie in $R$, a region of $D_2$ which
neighbours $v$, $\Omega_2^\uparrow$ moves may now be applied in turn
 to push the `nested tails' which have just been added over the two
edges of $R$ which are adjacent to $v$, as shown in \fullref{F:6}.

\begin{figure}[ht!] \centering
\includegraphics{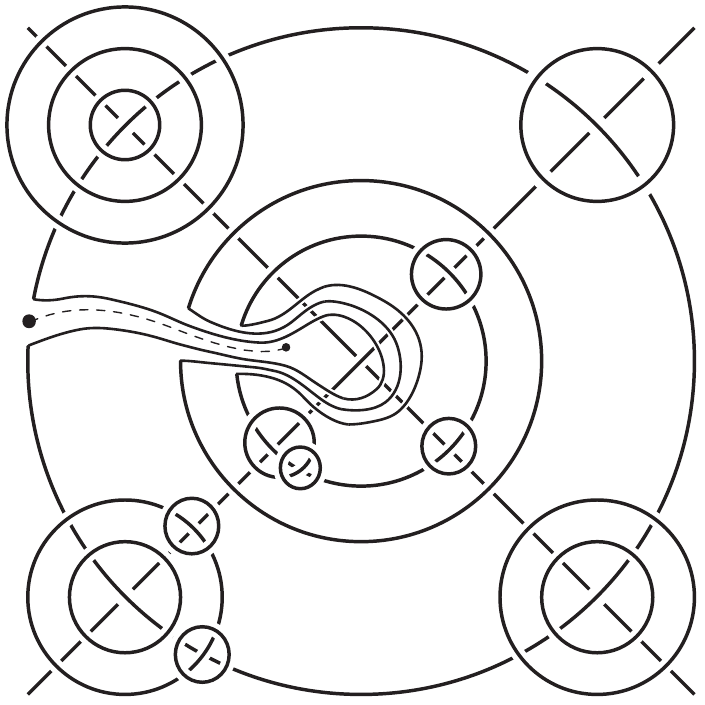}
\caption{Perform these $\Omega_3$ moves at each vertex}\label{F:7}
\end{figure}

Note that the procedure undertaken so far takes place inside $C$,
the outermost circle part associated with $v$, but outside any
circle parts associated to other vertices inside and on $C$. This is
good news since it means that we may repeat this operation at all
vertices with at least one circle part associated. After doing this,
we are done with $\Omega_2^\uparrow$ moves. We now form $D_3$ by
using $\Omega_3$ moves to push all the circle parts associated to a
particular vertex across that vertex, as in \fullref{F:7}, again
observing that we may do this on each collection of circle parts
independently.

\begin{figure}[ht!] \centering
\includegraphics{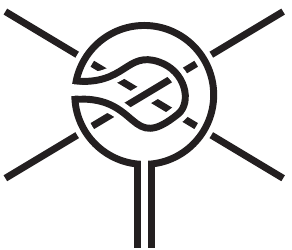}
\caption{Circle part of the last lollipop}\label{F:8}
\end{figure}

It is now time to show that $D_3$ may be obtained from $D_1$ by
means of a sequence of $\Omega_2^\uparrow$ moves. Let us go back to
the sequence $D_1 = E_1 , \ldots , E_n = D_2$ where  $E_i
\rightsquigarrow E_{i+1}$ or $E_i \circleright E_{i+1}$ for $i \in
\{1 , \ldots , n-1 \}$. Consider the part of $D_2$ which was added
in the final step. If this was a tail, then it still is in $D_3$ and
it may be removed by $\Omega_2^\downarrow$ moves. If it was a
lollipop then it may also now be removed by $\Omega_2^\downarrow$
moves since the circle part is now as shown in \fullref{F:8}.

After removing the last tail or lollipop from $D_3$ we may now
remove the second last in the same way. Repeating this process we
will eventually reach $D_1$ by means of $\Omega_2^\downarrow$ moves.
\endproof

\medskip
\textbf{\fullref{thm:1}}\qua {\sl Given two diagrams $D_1$ and $D_2$ for a
link $L$, $D_1$ may be turned into $D_2$ by a sequence of
$\Omega_1^\uparrow$ moves, followed by a sequence of
$\Omega_2^\uparrow$ moves, followed by a sequence of $\Omega_3$
moves, followed by sequence of $\Omega_2^\downarrow$ moves.

Furthermore, if $D_1$ and $D_2$ are diagrams of a link where the
winding number and framing of each component is the same in each
diagram, then $D_1$ may be turned into $D_2$ by a sequence of
$\Omega_2^\uparrow$ moves, followed by a sequence of $\Omega_3$
moves, followed by a sequence of $\Omega_2^\downarrow$ moves.}

\proof By \fullref{prop:1} it is enough to prove the second
part of the theorem. Let $D_1$ and $D_2$ be diagrams of a link where
the winding number and framing of each component is the same in each
diagram. Bruce Trace proved in \cite{Tr} that any two knot diagrams
with the same winding number and framing may be turned from one into
another by means of $\Omega_2$ and $\Omega_3$ moves. In fact his
result may be readily generalised to link diagrams with the same
hypotheses as we have made about $D_1$ and $D_2$. All one needs to
do is to apply the method used in \cite{Tr} to each component of the
link.

We shall thus proceed by induction on $M(D_1,D_2)$, the minimum
number of Reidemeister moves required to turn $D_1$ into $D_2$ with
only $\Omega_2$ and $\Omega_3$ moves. The claim clearly holds for
$M(D_1,D_2)=1$. Let $D_1=I_1,\ldots{},I_m=D_2$ be a sequence of link
diagrams arising from a minimal length sequence of $\Omega_2$ and
$\Omega_3$ moves connecting $D_1$ and $D_2$. Then
$M(D_1,D_2)=M(I_2,D_2)+1$. By the inductive hypothesis, $I_2$ may be
turned into $D_2$ by a sequence of $\Omega_2^\uparrow$ moves,
followed by a sequence of $\Omega_3$ moves, followed by a sequence
of $\Omega_2^\downarrow$ moves.

Let $$I_2=E_1,\ldots,E_{n}=D_2$$ be a sequence of diagrams arising
from such a sequence of Reidemeister moves. Let $E_p$, $E_q$
$(1\leq{}p\leq{}q\leq{}n)$ be such that
$$E_1 \stackrel{\Omega_2^\uparrow}{\longrightarrow} , \ldots , \stackrel{\Omega_2^\uparrow}{\longrightarrow} E_p  \stackrel{\Omega_3}{\longrightarrow} , \ldots ,
\stackrel{\Omega_3}{\longrightarrow} E_q
\stackrel{\Omega_2^\downarrow}{\longrightarrow} , \ldots ,
\stackrel{\Omega_2^\downarrow}{\longrightarrow} E_{n}.$$ If the move
from $I_1$ to $I_2$ is of type $\Omega_2^\uparrow$ then there is
nothing to prove. The remaining cases to consider are if this move
is of type $\Omega_3$ or of type $\Omega_2^\downarrow$. In the
former case apply \fullref{lem:1} to this move and the \fullref{cor:3}
to the sequence of $\Omega_2^\uparrow$ moves that follow it to
obtain a diagram $E_p'$ as shown:
$$\labellist\small
\pinlabel* {seq$(\Omega_2^\uparrow)$} [rb] <0pt,0pt> at 18 32
\pinlabel* {seq$(\Omega_2^\uparrow)$} [rb] <0pt,0pt> at 69 82
\pinlabel* {seq$(\Omega_2^\uparrow)$} [tl] <0pt,0pt> at 132 33
\pinlabel {$\Omega_3$} [t] <0pt,0pt> at 36 15
\pinlabel {$E_p'$} <0pt,0pt> at 97 108
\pinlabel {$E_p$} <0pt,0pt> at 150 56
\pinlabel {$E_1'$} <-1pt,1pt> at 43 53
\pinlabel {$I_2=E_1$} <0pt,0pt> at 96 15
\pinlabel {$I_1$} <0pt,0pt> at 2 15
\endlabellist
\includegraphics{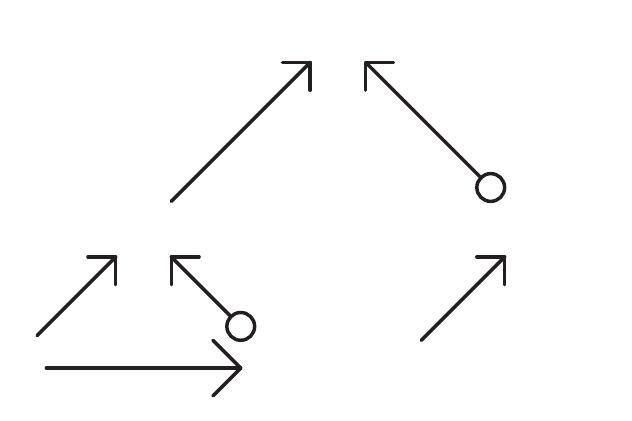}$$
If the move from $I_1$ to $I_2$ is of type $\Omega_2^\downarrow$
then $I_2 \rightsquigarrow I_1$. Thus we may apply the \fullref{cor:2}
to obtain $E_p'$ as shown:
$$\labellist\small
\pinlabel* {seq$(\Omega_2^\uparrow)$} [rb] <0pt,0pt> at 26 69
\pinlabel* {seq$(\Omega_2^\uparrow)$} [tl] <0pt,0pt> at 91 22
\pinlabel {$E_p'$} <0pt,0pt> at 57 98
\pinlabel {$E_p$} <0pt,0pt> at 112 41
\pinlabel {$I_2=E_1$} <3pt,2pt> at 50 0
\pinlabel {$I_1$} <0pt,0pt> at 0 46
\pinlabel \scalebox{2.7}{\rotatebox{135}{$\rightsquigarrow$}}  
<-3pt,-3pt> at 85 77
\pinlabel \scalebox{2.7}{\rotatebox{135}{$\rightsquigarrow$}}  
<6pt,6pt> at 14 17
\endlabellist
\includegraphics{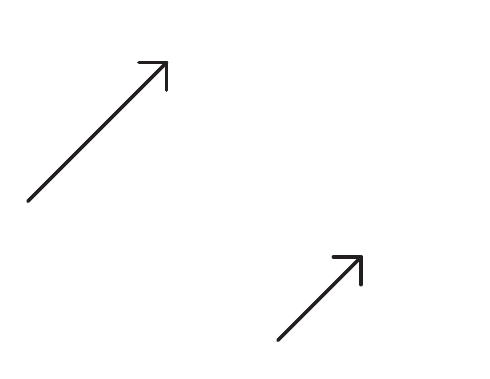}$$
Thus in either case we may perform a sequence of $\Omega_2^\uparrow$
moves on $I_1$ to obtain a diagram $E_p'$ such that $E_p
\circleright E_p'$ or $E_p \rightsquigarrow E_p'$. Now, $E_p$ and
$E_q$ are joined by a sequence of $\Omega_3$ moves. Applying \fullref{lem:1}
to each of these we obtain the following:
$$\labellist\small
\pinlabel* {seq$(\Omega_2^\uparrow)$} [rb] <0pt,0pt> at 20 35
\pinlabel* {seq$(\Omega_2^\uparrow)$} [rb] <0pt,0pt> at 100 35
\pinlabel* {seq$(\Omega_2^\uparrow)$} [rb] <0pt,0pt> at 227 35
\pinlabel {$\Omega_3$} [t] <0pt,0pt> at 36 14
\pinlabel {$\Omega_3$} [t] <0pt,0pt> at 117 14
\pinlabel {$\Omega_3$} [t] <0pt,0pt> at 247 14
\pinlabel {$E_p$} <0pt,0pt> at 2 14
\pinlabel {$E_{p+1}$} <1pt,0pt> at 79 14
\pinlabel {$E_{p+1}\ldots E_{q-1}$} <0pt,0pt> at 183 14
\pinlabel {$E_q$} <0pt,0pt> at 290 14
\endlabellist
\includegraphics{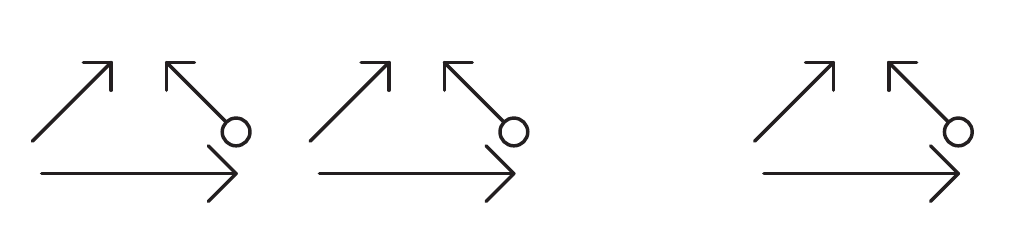}$$
The stage is now set to apply the \fullref{cor:3} several times
to obtain a diagram $E_q'$ as shown in the next diagram:

$$\labellist\small
\pinlabel* {seq$(\Omega_2^\uparrow)$} [rb] <2pt,1pt> at 18 24
\pinlabel* {seq$(\Omega_2^\uparrow)$} [rb] <2pt,1pt> at 57 64
\pinlabel* {seq$(\Omega_2^\uparrow)$} [rb] <2pt,1pt> at 98 24
\pinlabel* {seq$(\Omega_2^\uparrow)$} [rb] <2pt,1pt> at 139 145
\pinlabel* {seq$(\Omega_2^\uparrow)$} [rb] <2pt,1pt> at 179 105
\pinlabel* {seq$(\Omega_2^\uparrow)$} [rb] <2pt,1pt> at 220 66
\pinlabel* {seq$(\Omega_2^\uparrow)$} [rb] <2pt,1pt> at 258 26
\pinlabel {$E_p$} <0pt,0pt> at 1 6
\pinlabel {$E_{p+1}$} <2pt,0pt> at 78 6
\pinlabel {$E_{p+2}$} <2pt,0pt> at 159 6
\pinlabel {$E_q$} <0pt,0pt> at 321 6
\pinlabel {$E_{p+1}'$} <1pt,2pt> at 40 44
\pinlabel {$E_{p+2}'$} <0pt,0pt> at 80 86
\pinlabel {$E_q'$} <0pt,0pt> at 160 166
\pinlabel {\rotatebox{135}{$\cdots$}} <0pt,0pt> at 220 108
\pinlabel {\rotatebox{135}{$\cdots$}} <0pt,0pt> at 179 65
\pinlabel {\rotatebox{45}{$\cdots$}} <0pt,0pt> at 99 109
\pinlabel {\rotatebox{45}{$\cdots$}} <0pt,0pt> at 143 64
\pinlabel {\rotatebox{45}{$\cdots$}} <0pt,0pt> at 178 27
\pinlabel {$\cdots$} <0pt,0pt> at 200 6
\endlabellist
\includegraphics{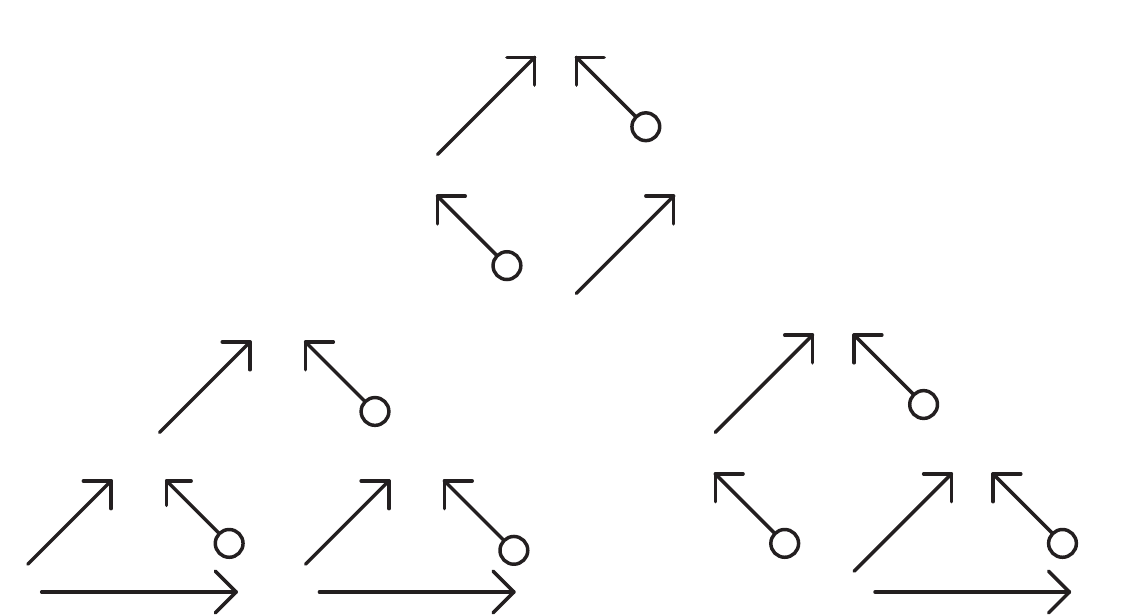}$$
Now, $E_p'$ is a diagram with either  $E_p \circleright E_p'$ or
$E_p \rightsquigarrow E_p'$. If we apply the \fullref{cor:2} or
the \fullref{cor:3} accordingly, then we get a diagram $E_q''$
as shown:
$$\labellist\small
\pinlabel {seq$(\Omega_2^\uparrow)$} [b] <0pt,0pt> at 31 74
\pinlabel {seq$(\Omega_2^\uparrow)$} [b] <0pt,0pt> at 193 74
\pinlabel* {seq$(\Omega_2^\uparrow)$} [t] <0pt,0pt> at 31 37
\pinlabel* {seq$(\Omega_2^\uparrow)$} [t] <0pt,0pt> at 193 37
\pinlabel {$E_p$} <0pt,0pt> at 2 37
\pinlabel {$E_p$} <0pt,0pt> at 160 35
\pinlabel {$E_p'$} <0pt,0pt> at 2 74
\pinlabel {$E_p'$} <0pt,0pt> at 160 74
\pinlabel {$E_q'$} <0pt,0pt> at 71 37
\pinlabel {$E_q'$} <0pt,0pt> at 238 35
\pinlabel {$E_q''$} <0pt,0pt> at 71 74
\pinlabel {$E_q''$} <0pt,0pt> at 238 74
\pinlabel \scalebox{2}{\rotatebox{90}{$\rightsquigarrow$}}  
<0pt,0pt> at 2 55
\pinlabel \scalebox{2}{\rotatebox{90}{$\rightsquigarrow$}}  
<0pt,0pt> at 71 55
\pinlabel {$I_1\stackrel{\Omega_2^\uparrow}{\longrightarrow}I_2$} at 35 13
\pinlabel {$I_1\stackrel{\Omega_3}{\longrightarrow}I_2$} at 197 13
\endlabellist
\includegraphics{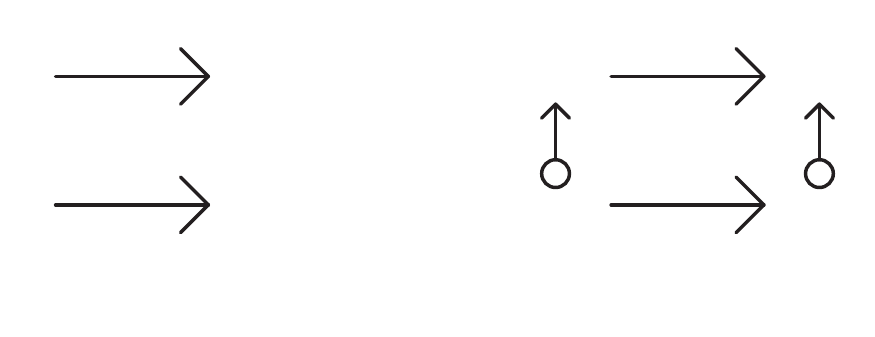}$$
Thus we have formed a diagram $E_q''$ such that:\begin{enumerate}
    \item $E_q''$ is obtained from $D_1$ by means of a sequence of
    $\Omega_2^\uparrow$ moves.
    \item $E_q''$ is obtained from $E_q$ by the addition of a
    sequence of tails and lollipops.
\end{enumerate}
We complete the proof by applying \fullref{prop:2} to $E_q''$ and
$E_q$. \endproof

We conclude this paper by noting that although the last theorem was
proved by induction, we could have taken any sequence of $\Omega_2$
and $\Omega_3$ moves as our ingredients. In this way, one could
obtain a (large) upper bound on the number of sorted moves required
to pass from one diagram to the other in terms of the minimum number
of unsorted moves.

\bibliographystyle{gtart}
\bibliography{link}

\end{document}